\documentclass [11pt]{amsart}
\usepackage{amscd}
\usepackage{amsxtra}
\usepackage{epsfig}
\usepackage{graphicx,color}
\usepackage{amssymb}
\usepackage[all]{xy}
\newtheorem{theorem}{\sc Theorem}
\newtheorem{proposition}{\sc Proposition}
\newtheorem*{definition}{\sc Definition}
 \makeatletter
    
    \@addtoreset{equation}{section}
  \makeatother  
\font\tmsb=msbm10 at12pt
\font\smsb=msbm7
\font\ssmsb=msbm5
\newfam\msbfam
\textfont\msbfam=\tmsb
\scriptfont\msbfam=\smsb
\scriptscriptfont\msbfam=\ssmsb

\def \RM{\mathbb {R}}%        corps des reels
\def \CM{\mathbb{C}}%        nombres complexes
\def \Kt {\mathcal{K}}
\def \HM {\mathbb{H}}

\def \Ct {\mathcal{C}}
\def \Et {\mathcal{E}}
\def \Ht {\mathcal{H}}

\def \Hom {\mathcal{H}om}
\def \d{\partial}%derivee partielle

\def\dt{\delta} 
\def\a{\alpha}

\def\p{\varphi}

\def\l{\lambda}
\def\L{\Lambda}
\def\lb{\left\{}
\def\rb{\right\}}

\def \t{\tilde}

\def \to{\longrightarrow} 
\def \w{\wedge}

\def\tensor{\hat{\otimes}}

\newcommand{\M}{{\mathcal M}}

\newcommand{\It}{{\mathcal I}}

\newcommand{\OM}{{\mathcal O}}

\newcommand{\Lt}{{\mathcal L}}

\title [A rigidity theorem for Lagrangian deformations]
{\sl  A rigidity theorem for Lagrangian deformations}
\author{\sl Mauricio D. Garay  }
\date{Novenber 2003}
\address{Fachbereich Mathematik 17\\
Staudingerweg 9, Johannes Gutenberg-Universit\"at\\
55099 Mainz, Germany.}
\email{garay@mathematik.uni-mainz.de}
\begin{document}
\parindent=0
cm
%%%%%%%%%%%%%%%%%%%%%%%%%%%%%%%%%%%%%%%%%%%%%%%%%%%%%%%%%
%%%%%%%%%%%%%%%%%% ABSTRACT%%%%%%%%%%%%%%%%%%%%%%%%%%%%%%
%%%%%%%%%%%%%%%%%%%%%%%%%%%%%%%%%%%%%%%%%%%%%%%%%%%%%%%%%

\begin{abstract}
We consider deformations of singular Lagrangian varieties in symplectic
spaces. We prove that a Lagrangian deformation of a Lagrangian complete intersection is analytically rigid provided that this is the case infinitesimally.
This result solves a problem posed by Colin de Verdi\`ere concerning Lagrangian curves. Finally, we prove the coherence of the direct image sheaves of relative infinitesimal Lagrangian deformations.

\end{abstract}
\maketitle
\section*{Introduction}
The study of deformations of quantum integrable systems in the semi-classical
limit is governed by the deformation of their characteristic varieties. 
Namely, to the symplectic equivalence of characteristic varieties 
corresponds the equivalence of the ideal generated by the 
corresponding micro-differential operators up to right multiplication 
by a micro-elliptic operator (see e.g. \cite{Colin_Parisse}).
For quantum integrable systems, the characteristic variety is a Lagrangian
variety which is as a general rule a singular complete intersection with non-isolated singularities.
The singularities of this variety reflect the complexity of the quantum 
system, for instance the singular Bohr-Sommerfeld rules or 
the asymptotic behaviour of the spectrum (see e.g. \cite{Parisse},  \cite{Sjostrand}).\\
%Singularities of Lagrangian varieties appear naturally in the  
%applications of symplectic geometry
%(\cite{swallow},\cite{Givental},\cite{Ishikawa},\cite{Zakalyukin}).
%In recent years, the study of semi-classical deformations of
%pseudo-differential operators attracted new attention  on them
%(\cite{Parisse}).
Following a previous work of Pham (\cite{Pham}, see also \cite{Kostov}),
Colin de Verdi\`ere proved a microlocal versal deformation theorem for
a particular class of one-dimensional microdifferential equations
(\cite{Colin}).
Using the above correspondance, the theorem was stated as
a Lagrangian versal deformation theorem for monomial deformations of quasi-homogeneous Lagrangian curve germs. Then, Colin de Verdi\`ere posed the problem of
the existence of a Lagrangian versal deformation theorem for general Lagrangian curve germs.
In this paper, we prove a rigidity theorem for
Lagrangian varieties which solves this problem positively when applied to the
case of Lagrangian curves.\\
%Our results are obtained by investigating the relative 
%infinitesimal Lagrangian deformations and by making a correspondence between
%coherence of direct images sheaves on and the
%rigidity, and versality theorems.
%This general fact is independent from the symplectic
%framework.
%It enables us to prove local and global
%rigidity theorems for infinitesimally versal Lagrangian deformations of
%Lagrangian complete intersection. Like in Colin de Verdi\`ere's paper,
%this results has an immediate interpretation in terms
%of semi-classical pseudo-differential operators.\\
This paper is organised as follows.\\ 
We start the first section by adapting to the relative case the deformation theory that was 
considered by Sevenheck and van Straten for 
Lagrangian singularities (\cite{VS}, \cite{Sev_these}). Following Sevenheck and van Straten, we introduce a complex of sheaves, called
the {\em Lagrange complex}, whose first 
cohomology sheaf computes the sheaf of local infinitesimal
Lagrangian deformation modulo trivial ones. Then we recall the main result of
\cite{mutau} namely the freness of the direct image sheaves of the relative Lagrange complex for infinitesimally Lagrangian versal deformations.\\
For smooth Lagrangian manifolds, the Lagrange complex is canonically isomorphic to the holomorphic de Rham complex.
For general Lagrangian varieties the relation between the Lagrange complex and the de Rham complex is more 
subtle and will not be investigated in this paper. \\

In the second section, we state the rigidity theorem and its main corollary which answers positively to the problem of Colin de Verdi\`ere.
This rigidity theorem is a weaker version of a conjectural Lagrangian versal deformation
theorem. We point out that we consider only a naive viewpoint of versality
namely versality over a smooth basis. We choose the name rigidity rather than
stability because we reserve the notion of stability to momentum mappings
(\cite{moment}).\\
In the third section, the proof of the rigidity theorem is given.\\
In the fourth section, we state and prove the coherence of the direct image sheaves of
the relative Lagrange complex. The proof is based 
on a practical corollary of the Kiehl-Verdier theorem. Once this corollary is established, we just adapt  to the relative case the proof given by Sevenheck and van Straten for the finiteness of the cohomology sheaves of the Lagrange complex (\cite{VS}). This coherence result is used in the proof of the rigidity theoremand in the proof of the main theorem of \cite{mutau}.\\
 In the appendix, we prove a corollary to the
Kiehl-Verdier theorem which was used in section 4 to establish the coherence
theorem.\\
The results of this paper can be easily adapted for real analytic  Lagrangian varieties.
The formal aspect of the theory is treated in Sevenheck's thesis (\cite{Sev_these}).
%The formal theory is more 
%simple and as Sevenheck proved in his thesis the Lagrangian versal 
%deformation theorem holds in the formal category. In fact,
%Schlessinger's theory can be entirely reformulated in the Lagrangian
%context (see \cite{Sev_these} for details).A preliminary version
%of this paper appeared in \cite{hamilton}.\\
\ \\

{\em Acknowledgements.}
% This work has been realised at the Department 
%of Algebraic Geometry in Mainz 
%University and at the Laboratoire de G\'eom\'etrie et Dynamique at the
%Institut de Math\'ematiques de Jussieu. The 
%author would like to thank both of these institutions for their hospitality.
It is a pleasure to thank D. van Straten for numerous 
discussions from which this work originated and C. Sevenheck for his explanations
on singular Lagrangian varieties. 
The author also acknowledge M. Chaperon, P. Dingoyan, H. Eliasson
E. Ferrand, B. Teissier, C. Voisin. Last
but not least, the author thanks Y. Colin de Verdi\`ere whose work
and questions motivated the present research.
This work was supported by the Deutsche 
Forschungsgemeinschaft (Project no. SPP 1094).

%%%%%%%%%%%%%%%%%%%%%%%%%%%%%%%%%%%%%%%%%%%%%%%%%%%%%%%%%%%%%
\section{Lagrangian deformations}
%About 10 years ago, van Straten introduced a complex of sheaves whose first 
%cohomology parameterise the infinitesimal deformations of 
%Lagrangian varieties.
%With the help of this complex, he proved the
% rigidity of the  open Whitney umbrella in $ \CM^4 $
% (i.e.: the conormal bundle to the curve $ \lb y^2=x^3 \rb $) as a variety.
% For Lagrangian immersions this
% theorem is also true
The complex of infinitesimal Lagrangian deformations
is a particular case of the complex which 
computes the deformations of Lie algebroids (see e.g. \cite{Mackenzie}).
It was introduced for the study of Lagrangian singularities in \cite{VS}.
For a detailed exposition of the deformation theory of Lagrangian 
singularities, we refer to Sevenheck's thesis (\cite{Sev_these}).
%%%%%%%%%%%%%%%%%%%%%%%%%%%%%%%%%%%%%%%%%%%%%%%%
\subsection{Deformations of real compact Lagrangian manifolds}
\label{SS::compact}
To understand the construction of the Lagrange complex, it is 
useful to investigate the 
deformations of a smooth compact Lagrangian submanifold
$ L \subset \RM^{2n} $.
For such a manifold, the Darboux-Weinstein theorem asserts that 
there exists a symplectomorphism
$$ \p:\RM^{2n} \to T^*L $$ which maps a tubular neighbourhood
of $ L \subset \RM^{2n} $ to a tubular neighbourhood of the zero 
section in the cotangent bundle $ T^*L $ to $ L $ 
(see e.g. \cite{Weinstein}).\\
Via the map $ \p $, a small one-parameter deformation $(L_t)$ of $ L $ is 
mapped to the family of graphs of a one parameter family of maps
$$ \a_t:L \to T^*L $$ that is, to a family of differential 
one forms.\\ It is readily verified that $L_t$ is Lagrangian if 
the one-form $ \a_t $ is closed and that
$ L_t $ is Hamiltonian isotopic to $ L_0=L $ if $ \a_t $ is exact.\\
Consequently, the space of  
deformations of $ L $ over some base $ \L $ modulo Hamiltonian isotopies
is parameterised by the maps from $\L$ to the first de Rham cohomology group
$H^1(L,\RM)$ of $ L $ sending $ 0 
\in \L $ to $ 0 \in H^1(L,\RM) $.
For the global study of deformations of complex Lagrangian manifolds, we refer
to \cite{Voisin}.
%For (smooth) complex Lagrangian manifolds, Voisin constructed a deformation
%theory using the holomorphic de Rham complex (\cite{Voisin}).
%It was proved in that paper that for smooth complex Lagrangian manifolds
%there exists a smooth Lagrangian versal
%deformation space. 
%For singular Lagrangian varieties, the relationship between the de Rham
% complex and the Lagrangian deformations is more subtle.
% and uses the 
% Lagrangian Gauss-Manin connection (\cite{Gauss}). It is precisely
% because of these difficulties that van Straten proposed to study another
% complex directly related to the Lagrangian deformations. 
% \begin{figure}[ht]
% \begin{center}
% \input{cylindre.pstex_t}
% \end{center}
% \caption{\it Here the Lagrangian manifold $L=L_0$ is the zero section of $ T^*S^1$.
% The space of infinitesimal deformations modulo Hamiltonian isotopies is
% generated by the restriction to $t=0$ of $\d_t \p_t$ where
% $\p_t$ is the translation by $t$ along the fibres of $ T^*S^1 \to S^1$. The
% translated Lagrangian manifold $L_t$ is the graph of the one-form $td \theta$
% where $\theta \in S^1$ is a coordinate on the zero section of $ T^*S^1$. The
% differential one-form $d \theta$ generates the first de Rham cohomology
% $H^1(L)$ of $L$.}
% \label{F::1}
% \end{figure}
%%%%%%%%%%%%%%%%%%%%%%%%%%%%%%%%%%%
\subsection{The relative Lagrange complex}
\label{SS::complex}
%We define the relative version of the construction of the complex of 
%infinitesimal Lagrangian deformations (\cite{VS}).\\
We consider a complex manifold
$M$ of dimension $2n$ together with a holomorphic symplectic 
two form $ \omega \in \Omega^2_{M}(M) $.\\
The Poisson bracket $ \lb f,g \rb $ of two 
holomorphic functions $ f,g \in \OM_M(U) $ is defined by the formula
$$\lb f,g \rb \omega^n=df \w dg \w \omega^{n-1}.$$
Recall that a {\em Lagrangian submanifold} of $M $ is an $ n 
$-dimensional holomorphic manifold on which the symplectic form 
vanishes. 
By Lagrangian variety $ L \subset M $, we mean a
reduced complex space of pure dimension $ n $ defined by an ideal sheaf $ \It_L $ such that 
which is closed under the Poisson bracket.
This means that
$$\lb f,g \rb \in \It_L(U)  $$
whenever $ f,g \in \It_L(U) $.\\

We consider now the situation with parameters.
Let $ \L $ be a complex space (the parameter space) with a 
marked point, denoted by $ 0 $, on it.  Unless specific mention, we assume that
$\L$ is a smooth complex manifold. The 
Poisson bracket on $ M $ lifts to an $\OM_\L$-linear  Poisson bracket on
$ M \times \L $.  We shall say that a variety
$$ L \subset (M \times \L)$$
is a Lagrange variety if its projects one-to-one to a Lagrange variety 
of $ M $. 
Let $ Z \subset (M \times \L) $ be a reduced complex subspace with ideal 
sheaf $ \It $. A map
$$ \p:Z \to \L $$
is called a {\em Lagrangian deformation} of its $ 0 $-fibre if
$ \p $ is a flat deformation and if the fibres of $ \p $ are reduced Lagrangian
varieties.\\
For a given $ f=f_1 \wedge \dots \wedge f_k \in \bigwedge^k \It/\It^2 $,
we denote by $ f^j \in \bigwedge^{k-1} \It/\It^2$  the element  
$$ f^j=(-1)^j f_1 \wedge \dots \wedge {\hat f_j} \wedge \dots \wedge
f_k . $$
The element $ (f^i)^j \in \bigwedge^{k-2} \It/\It^2 $ is denoted by $ f^{i,j} $.

\begin{definition}{\rm The {\em relative Lagrange complex}, denoted $ 
(\Ct_{Z/\L}^\cdot,\dt) $, of the Lagrangian deformation
$ \p:Z \to \L $ is the complex of sheaves on $Z$ defined by 
$$ \Ct_{Z/\L}^k={\Hom}_{\OM_Z}(\bigwedge^k \It/\It^2,\OM_Z) $$
and the $ k^{th} $ differential $ \dt:\Ct_{Z/\L}^k \to \Ct_{Z/\L}^{k+1} $ is given by
$$ \dt[\p](f)=\sum_{1 \leq i \leq k}\lb f_i,\p(f^i) \rb-\sum_{1 \leq 
i<j \leq k}\p(\lb f_i,f_j \rb \w f^{i,j}). $$}
\end{definition}
{\em Notation.}{ If $ L $ is a Lagrange variety then the {\em Lagrange complex}
of $ L $, denoted by $ \Ct_L^\cdot $ (which is defined in \cite{VS}) is the 
relative Lagrange complex of the constant deformation
$$\p:L \to\lb 0 \rb.  $$}
% {\em Remark.}{ The first cohomology vector space of the complex $ \Ct_L^\cdot $
% describes the infinitesimal deformations of $ L $. Namely, the normal 
% sheaf $ \Hom_{\OM_L}(\Lt,\OM_L) $ is usually interpreted as the space of
% infinitesimal deformation of a variety (see e.g. \cite{Pfister}).
% The condition $ \dt \p=0 $ ensures that the deformation is 
% infinitesimally Lagrangian and $ \p=\dt(h) $ that it can be trivialised
% by the symplectomorphism which is the time one map of the flow of $ h $. 
%Nevertheless, it seems that unlike in the usual deformation theory,
%the higher order cohomology space should not be considered as reflecting
%the obstructions to extend deformations but rather the topology of a
%smoothing of the 
%Lagrangian variety, when it exists (\cite{mutau}).
%}

%%%%%%%%%%%%%%%%%%%%%%%%%%%%%%%%%%%%%%%%%%%
%%%%%%%%%%%%%%%%%%%%%%%%%%%%%%%%%%%%%%%%%%%
\subsection{Equivalence of Lagrangian deformations}
\label{SS::ci}
%Recall that a Hamiltonian isotopy
%$$g_t:M \to M$$
%is an isotopy obtained by integrating the flow of the Hamiltonian vector field
%of some time-dependent holomorphic function
%$$H:M \times[0,1] \to \CM.$$
\begin{definition}{\rm A Lagrangian deformation-germ
$\p':(Z',0) \to (\L',0)$ is 
called
{\em $ \Lt $-induced} (resp. {\em $ \Lt $-equivalent}) from (resp. to)
$ \p:(Z,0) \to (\L,0) $ if there exists a Poisson mapping
(resp. a biholomorphic Poisson mapping)
$ g $ which makes the following diagram commute
$$ \xymatrix{(Z',0) \ar[r]^g \ar[d]^{\p'}& (Z ,0)\ar[d]^\p \\
              (\L',0) \ar[r] & (\L,0)} $$}
\end{definition}
%In the sequel, unless specific mention, we consider only deformation 
%over smooth basis.
\begin{definition}
\label{D::rigid}
{\rm A Lagrangian deformation-germ $ \p:(Z,0) \to (\L,0) $ is called
{\em rigid} if any deformation-germ $ \p':(Z',0) \to (\L \times T,0) $ such
that the restriction of $ \p'$ above $\L \times \lb 0 \rb$ equals 
$\p $, is induced from $ \p $.}
\end{definition}
A deformation
$$ \p:(Z,0) \to (\L,0),\ \ \ (\L,0) \approx (\CM^k,0) $$
of a Lagrangian variety $ L $ is called {\em $ \Lt $-versal},
if any deformation of $ L $ over any smooth basis is induced from $ \p $.\\
%In the sequel, we shall consider the germ of a deformation along
%a fibre $ L $ at a point, say $ 0 \in \L $
%$$ \p:Z \to (\L,0) .$$
%This means that we localise $ \OM_Z $ at the ideal sheaf $ \It_{L} $. 
% Let us spell out this notion.\\
%Two deformations $\p_1,\p_2$ have the same germ at $0 \in \L$ if there exist
%a neighbourhood $U \subset (Z_1 \cap Z_2)$, such that the restrictions of the
%deformations to $U$ are equal. The equivalence class is denoted by
%$$ \p:Z \to (\L,0).$$
%In the language of schemes $ \p$ is the germ of $\p_1$ at the ''point''
%$\p_1^{-1}(0)$.\\
After a choice of local coordinates in $ M \times \L $ at a point $ x $,
the ideal of $Z$ at $x$ is generated by holomorphic function-germs
$$F_1,\dots,F_p:(\CM^k \times \CM^{2n},0) \to (\CM,0).$$
We say that $F=(F_1,\dots,F_p)$ is the deformation of the
germ of $L$ at $p$, where $L \times \lb 0 \rb$ is the $0$-fibre of
$\p$.\\

%%%%%%%%%%%%%%%%%%%%%%%%%%%%%%%%%%%%

%%%%%%%%%%%%%%%%%%%%%%%%%%%%%%%%%%%%%%%%%%%%%%%%%%%%%%%%%%%%%%%%
\subsection{The Lagrangian Kodaira-Spencer map}
Consider a deformation
$$F:(\CM^k \times \CM^{2n},0) \to (\CM^n,0),\ \ \ (\l,x) \mapsto
(F_1(\l,x),\dots,F_n(\l,x))$$
of a Lagrangian complete intersection germ $L$.
As $L$ is a complete intersection, we get identifications
$$\Ct^k_{Z/\L} \approx \bigwedge^k \OM_{Z/\L}^n$$
given by the choice of the $F_i$'s.\\
We define a {\em Kodaira Spencer map} by 
$$\begin{matrix}\theta_F:& T_0\CM^k &\to& H^1(\Ct^{\cdot}_{Z/\L,0})\\ 
&v &\mapsto &[(DF.v )]
\end{matrix}$$
For instance 
$$ \theta_F(\d_{\l_i})=[(\d_{\l_i}F_1,\dots,\d_{\l_i}F_n )] .$$
This Kodaira-Spencer map is well-defined,
indeed differentiating the following equality along $v$
$$\lb F_j,F_k \rb=\sum_{i=1}^n a_i F_i$$
we get that
$$\lb d F_j.v,F_k \rb+\lb F_j,d F_k.v \rb=\sum_{i=1}^n (d a_i.v) F_i+
\sum_{i=1}^n a_i (d F_i.v)$$
which implies that the mapping $DF.v $ defines a coboundary in
$\Ct^\cdot_{Z/\L,0}$.
\begin{definition}{\rm The map $\theta_F:T_0\CM^k \to 
H^1(\Ct^{\cdot}_{Z/\L,0})$
defined above is called the  {\em relative Lagrangian Kodaira-Spencer map}
of $F$.}
\end{definition}
By restriction to $\l=0$, we define an {\em absolute Kodaira-Spencer mapping}: 
$$\begin{matrix}\bar \theta_F:& T_0\CM^k &\to& H^1(\Ct^{\cdot}_{L,0}),\\ 
&v &\mapsto &[(DF.v )_{\mid \l=0}]
\end{matrix}$$
and a commutative diagram
$$\xymatrix{
& H^1 (\Ct^{\cdot}_{Z/\L,0}) \ar[d]^{r} \\
T_0\L \ar[ru]^-{\theta_F}  \ar[r]^-{\bar \theta_F}  & H^1(\Ct^{\cdot}_{L,0})}$$
We denote by $ \Lt T^1_L $ the sheaf on $L$ of deformations of a 
Lagrangian variety $ L \subset M$ over the double 
point $ Spec(\CM[t]/t^2) $ modulo $ \Lt $-equivalence.
The proof of the following proposition is straightforward.
\begin{proposition}
[\cite{VS}]
\label{P::T1}
{The Lagrangian Kodaira-Spencer map gives an identification between the sheaf
$ \Lt T^1_L $ and the cohomology sheaf $ \Ht^1(\Ct^{\cdot}_L) $.}
\end{proposition}

\subsection{Pyramidal deformations}
\label{SS::pyramidal}
Following \cite{VS}, we introduce a class of singularity which plays 
the role of isolated singularities in symplectic geometry.\\
On the total space $ Z $ of the deformation
$$ \p:Z \to \L   $$
of a Lagrangian variety $ L \subset M $ a 
stratification is defined as follows:\\
Let $ f_1,\dots,f_k $ be generators
of the ideal of $ Z $ at a point $ x $. Denote by $ X_1,\dots,X_k $ be the 
Hamilton vector-fields of $ f_1,\dots,f_k $ and by $ l(x) $
the dimension of the vector spaces generated by the $ X_i $'s at $ x  $.
The strata $Z_j$ is defined by:
$$ Z_j=\lb x \in Z: l(x)=j \rb. $$
We have $Z=\bigcup_{j=0}^n Z_j$ where $2n=dim(M)$.\\
The following notion was introduced in \cite{VS} for Lagrangian 
varieties where it is called ''condition (P)''.
\begin{definition}
{\rm A Lagrangian deformation $ \p:Z \to \L $ is called {\em pyramidal}
if for any $ k $ the variety $ Z_k $ is of relative dimension at most $ k $.}
\end{definition}

 \begin{proposition}
\label{P::produit}
{The germ of a pyramidal deformation
 $ \p:Z \to \L $ at a point $ x \in Z_k$ is $\Lt$-equivalent
to a deformation-germ of the type
 $$\p':(Z' \times \CM^k,0) \to (\L,0),\ Z' \subset \CM^{2n-2k},\ 
 \CM^k \subset \CM^{2k}  $$
 which is constant on the second factor. Moreover, this decomposition
induces a quasi-isomorphism between the complexes $\Ct^\cdot_{Z/\L,x}$
and  $\Ct^\cdot_{Z'/\L,0}$.}
 \end{proposition}
The proof of this proposition is straightforward, it is based on a simple
symplectic reduction argument (in the absolute case see e.g. \cite{VS}).
\begin{theorem}[\cite{VS}]
\label{T::VS}{If $ (L,0) $ is the germ of a pyramidal 
Lagrangian variety then the cohomology space $ H^k(\Ct_{L,0}^\cdot) $ are 
finite dimensional vector spaces.}
\end{theorem}
{\em Conjecture.}{ If $\Ht^1(\Ct^\cdot_L)$ is finite and if $ L $ is a 
complete intersection then $L$ is pyramidal.}\\

Finally, we recall the main theorem of \cite{mutau}.
\begin{theorem}[\cite{mutau}]
\label{T::mutau}{Let $ \p:(Z,0) \to (\L,0) $
be the germ of a Lagrangian deformation of a pyramidal Lagrangian variety $ 
L $.
Assume that the Kodaira-Spencer map
$$ \bar \theta_F:T_0\L \to H^1(\Ct^\cdot_{L,0}) $$
is surjective. Then, the $\OM_{\L,0}$-module $ H^1(\Ct^\cdot_{Z/\L,0}) $ is free
and we have a canonical isomorphism 
$$ H^1(\Ct^\cdot_{Z/\L,0}) \approx  \OM_{\L,0} \otimes H^1(\Ct^\cdot_{L,0}).$$}
\end{theorem}
The proof of this theorem used the coherence theorem
(Theorem \ref{T::coherent},section \ref{S::coherent}). Therefore,
to avoid cross references it will not be used.

\section{The rigidity theorem}
%%%%%%%%%%%%%%
\subsection{Lagrangian versal deformation theorem of Colin de Verdi\`ere}
The following result is, according to Colin de Verdi\`ere, ``the main non-trivial result'' of \cite{Colin}.
 \begin{theorem}
[\cite{Colin}]
\label{T::Colin}
{A monomial Lagrangian deformation $ F:(\CM^k \times \CM^{2},0) \to (\CM,0) $
 of a quasi-homogeneous Lagrangian curve-germ $(L,0) \subset (\CM^{2},0)$ with an 
 isolated singular point is $ \Lt $-versal provided that the absolute Kodaira-Spencer map associated to  $F$ is surjective.}
\end{theorem}
{\em Remark (see e.g. \cite{Colin} or \cite{isochore} for details).}{ The surjectivity of the Kodaira-Spencer map is easily 
verified. Namely, denote by $ f $ a reduced equation of $ L $.
Then, as shown by Brieskorn (\cite{Br3}) there is a canonical identification
$$H^1(\Ct^{\cdot}_{L,0}) \approx \OM_{\CM^2,0}/Jf$$
induced by the identity mapping. Here $ Jf $ denotes the Jacobian ideal of $ f $.\\
Take a monomial basis of $\OM_{\CM^2,0}/Jf$ and lift it to
$$e_1,\dots,e_{\mu} \in \OM_{\CM^2,0}.$$
The quoted theorem of Colin de Verdi\`ere states that the deformation
$$f+\sum_{i=1}^\mu \l_i e_i$$
is $\Lt$-versal.}
%%%%%%%%%%%%%%%%%%%%%%%%%%%%%%%%%%%%%%%%%
\subsection{The rigidity theorem}
\label{SS::rigid}
\begin{theorem}
\label{T::versel}
{Let $ F:(\CM^k \times \CM^{2n},0) \to (\CM^n,0) $
be a Lagrangian deformation  of a  Lagrangian complete intersection-germ
$(L,0) \subset (\CM^{2n},0)$. Assume that
\begin{enumerate}
\item $ (L,0) $ is pyramidal,
\item the absolute Kodaira-Spencer map
$$\bar \theta_F:T_0\CM^k \to H^1(\Ct_{L,0})$$
associated to  $F$ is surjective.
\end{enumerate}
Then, the Lagrangian deformation $ F $ is rigid.}
\end{theorem}
{\em Example.}{ Consider the germ at the origin of the Lagrangian variety 
$$L= \lb (q,p) \in \CM^{2n}: q_1p_1=q_2p_2=\dots=q_np_n=0 \rb $$
and let $F=(F_1,\dots,F_n)$ be the deformation defined by
$$F_i=q_ip_i+\l_i.$$
 The Lagrangian Kodaira-Spencer maps  $\d_{\l_1},\dots,\d_{\l_n}$ to
the cohomology classes
$$[(1,0,\dots,0)],\dots,[(0,\dots,0,1)].  $$
A straightforward computation shows that 
they generate $H^1(\Ct^{\cdot}_{L,0}) $.\\
Thus, the Lagrangian deformation $F=f+(\l_1,\dots,\l_n)$
is rigid. This is a variant of the famous Vey-R\"ussmann theorem for separable
integrable systems (\cite{Russmann}, \cite{Vey}, \cite{Eliasson}). (In 
fact using the results of this paper together with that of \cite{moment}, 
a new proof of this theorem can be derived.)\\

The following result answers positively to Colin de Verdi\`ere's problem.
\begin{theorem}
\label{C::versel}
{A Lagrangian deformation $ F:(\CM^k \times \CM^{2},0) \to (\CM,0) $
 of a Lagrangian curve-germ $(L,0) \subset (\CM^{2},0)$ with an 
 isolated singular point is
$ \Lt $-versal provided that the absolute Kodaira-Spencer map associated to  $F$ is surjective.}
\end{theorem}
\begin{proof}
We follow standard arguments due to Martinet for the case of 
singularity theory of differentiable maps (\cite{Mar}).\\
Put $f=F(0,.)$ and let $ G $ be an arbitrary $ s $-parametric Lagrangian deformation of
$ f=F(0,.) $. We have to prove that $ G $ is $\Lt$-induced from $ F $.
To do this, we define the sum of $ F $ and $ G $ by the formula
$$ (F \oplus G)(\a,\l,x)=F(\l,x)+G(\a,x)-f .$$
The restriction of 
this deformation to $ \l=0 $ is equal to G. Consequently it is 
sufficient to prove that $ F \oplus G $ is $\Lt$-induced from $ F $.
Denote by $ F_j $ the restriction of $ F \oplus G $ to $ 
\a_j=\dots=\a_s=0 $. We have $ F_1=F $ and $ F_s=F \oplus G $.\\
The rigidity theorem implies that $ F_j $  is $\Lt$-induced 
from $ F_{j-1} $. By induction, we get that $ F_{j} $ is $\Lt$-induced 
from $ F_{j-k} $.
In particular, $ F_s= F \oplus G $ is $\Lt$-induced from $ F_1=F $.
\end{proof}
{\em Remark.} In case $n > 1$, we cannot use Martinet's argument
since $ F \oplus G$ is in general not a Lagrangian deformation.

%%%%%%%%%%%%%%%%%%%%%%%%%%%%%%%%%%%%%%%%%%%%%%%%%%%%%%%%%%%%
\section{Proof of the rigidity theorem.}
\label{S::versel}
%%%%%%%%%%%%%%%%%%%%%%
%%%%%%%%%%%%%%%%%%%%%%%%%%%%%%%%%%%%%%%%%%%%%%%5
%%%%%%%%%%%%%%%%%%%%%%%%%%%%%%%%%%%%%%%%%%%%%%%%%%%%%%%%%%%%%%%%%%
\subsection{Infinitesimal formulation of the problem}
\label{SS::homological}
Let
$$G:(\CM \times \CM^{k} \times \CM^{2n},0) \to (\CM^n,0),\ \ \ 
(t,\l,x) \mapsto G(t,\l,x)$$  be a deformation
of
$$F:(\CM^{k} \times \CM^{2n},0) \to (\CM^n,0),\ \ \ 
(\l,x) \mapsto F(\l,x).$$
The aim of this subsection is to prove the following assertion.\\
Assertion.{\em The deformation $G$ is $ \Lt $-induced from $F$ provided that there exist 
function-germs $ h \in \OM_{2n+k+1},b_1,\dots,b_k \in 
\OM_{k+1} $ and a matrix $ B \in gl(n,\OM_{2n+k+1})$ solving the
equation
\begin{equation}
\label{E::homologique}
 \lb h,G \rb+ B G +\sum_{i=1}^k b_i\d_{\l_i}G=-\d_t G .
\end{equation}}
  We search for a relative symplectomorphism-germ
$$ \p:(\CM \times \CM^k \times \CM^{2n} ,0) \to
(\CM^{2n},0),\ \ \ (t,\l,x) \to \p(t,\l,x)  $$ and
a matrix $ A \in GL(n,\OM_{2n+k+1})  $ such that the following 
equalities hold
\begin{equation}
\label{E::dbt1}
\lb 
\begin{matrix}F&=&A_{\tau}  ( G_\tau \circ \p_{\tau}),\\
(A_0,\p_0)&=&(I,Id).
\end{matrix}
\right.
\end{equation}
We have used the standard notations $I \in GL(n,\OM_{2n+k+1})$ for the
identity matrix, $Id$ for the identity mapping in $\CM^{2n}$, $ 
A_{\tau} $ for $ A(\tau,.,.) $ and so on.\\
Differentiating the first equation of the system (\ref{E::dbt1}) with 
respect to $ \tau $
 at $ \tau=t $, 
we get the equation
\begin{equation}
\label{E::db}
 A_t  \frac{d}{d\tau}_{\mid_{\tau=t}} (G_{t} \circ \p_\tau)+
 (\frac{d}{d\tau}_{\mid_{\tau=t}}  A_\tau )(G_{t} \circ 
\p_t)+A_{t} 
  (\frac{d}{d\tau}_{\mid_{\tau=t}}
G_\tau) \circ \p_{t}=0 .
\end{equation}
Define the time-dependent vector field germ $ v_{t} $ and the matrix 
$ B \in gl(n,\OM_{2n+k+1})  $ by the 
formulas
$$\lb \begin{matrix}
v_{t}(\p_{t}(\l,x))&= &\frac{d}{d\tau}_{\mid_{\tau=t}} \p_\tau(\l,x),\\
(A_t    B_t) \circ \p_t&=&\frac{d}{d\tau}_{\mid_{\tau=t}} A_\tau .
\end{matrix} \right.$$
 
Multiplying equation (\ref{E::db}) on the right by $ \p_t^{-1} $ and 
on the left by $ A_t^{-1} $, we get the equation
\begin{equation}
\label{E::homo}
 L_{v_t} G_{t}+ B_t  G_t+\d_t G_t=0 .
\end{equation}
Standard theorems on differential equations imply $ (\p_\tau,A_\tau) $ 
satisfying the system (\ref{E::dbt1}) can be found provided that 
there 
exist $ (v_t,B_t) $ satisfying equation (\ref{E::homo}).
Until here our arguments have been standard and hold for most of the
versal deformation theorems. We 
now come to the specificity of our situation.
Because the vector
field $v(t,.)=v_t$ comes from a relative symplectomorphism-germ, it is of 
the type
$$v=\sum_{i=1}^{2n} a_i\d_{x_i}+\sum_{i=1}^k b_i\d_{\l_i},\ 
a_i \in \OM_{2n+k+1},\ b_i \in \OM_{k+1}$$
where the vector field $w= \sum_{i=1}^{2n} a_i\d_{x_i} $ is the Hamiltonian
field of some  function-germ $ h \in \OM_{2n+k+1} $. 
Consequently, equation (\ref{E::homo}) can be written in the form
\begin{equation}
 \lb h,G \rb+ B  G +\sum_{i=1}^k b_i\d_{\l_i}G=-\d_t G .
\end{equation}
This proves our assertion.

%%%%%%%%%%%%%%%%%%%%%%%%%%%%%%%%%%%%%%%%%%%%%%%%%%%%%%%%%%%%%
\subsection{Solving the infinitesimal equation}
\label{SS::solve}
We interpret Equation (\ref{E::homologique}) in cohomological terms.\\
Let $ \p:Z \to \L $ and $ \p':Z' \to \L' $ be representatives for $ F $ and $ G $.
We have the following commutative diagram
$$\xymatrix{T_0\L' \ar[r]^-{\theta_G} \ar[d]^{\pi} & H^1 (\Ct^{\cdot}_{Z'/\L',0}) \ar[d]^{r_1} \\
T_0\L \ar[r]^-{\theta_F} \ar@/^/^{i}[u] \ar[dr]^{\bar \theta_F}& H^1 (\Ct^{\cdot}_{Z/\L,0}) \ar[d]^{r_2} \\
  & H^1(\Ct^{\cdot}_{L,0})}$$
The maps $\theta_G$ and $\theta_F$ are relative Kodaira-Spencer
maps, while $\bar \theta_F$ is the absolute one.
The maps $\pi,i$ are the canonical projection and injection induced by the
product structure $(\L',0) \approx (\L \times \CM,0)$. The map $r_1$ is the
restriction to $t=0$ and the map $r_2$ is the restriction to $\l=0$.\\
In this setting, Equation (\ref{E::homologique}) 
can be rewritten as
\begin{equation}
\label{E::fine}
 \sum_{i=1}^k b_i \theta_G(\d_{\l_i})=-\theta_G(\d_t), \ \ \ \  b_1,\dots,b_k 
\in \OM_{\L',0}.
\end{equation}
(We have used the notation $ \theta_G(\d_{\l_i}) $ rather than $ [\d_{\l_i}G] $
to underline in which space we take the cohomology class.)\\
Equation (\ref{E::fine}) can be solved provided that the cohomology class $\theta_G(\d_t)$ belongs to the $\OM_{\L',0}$-module
generated by the $\theta_G(\d_{\l_i})$'s.\\

As $\bar \theta_F$ factors through $\theta_F$, the restriction to zero map
$$ R:H^1 (\Ct^{\cdot}_{Z'/\L',0})/(\M_{\L',0} \otimes_{\OM_{\L',0}}H^1 (\Ct^{\cdot}_{Z'/\L',0}))  \to
H^1 (\Ct^{\cdot}_{L,0})$$
is surjective.\\
Assertion. {\em The map $R$ is an isomorphism.}\\
The assertion follows from Theorem \ref{T::mutau}
(section \ref{SS::pyramidal}) but in order to be
self-contained, we prove it.\\
Let us denote by
$$\p_p:Z_p \to \L_p,\ \ \ p=-1,0,\dots,k$$
the restriction of $\p'$ above 
$$\L_p=\lb (t,\l) \in \L':\l_0=\dots=\l_p=0 \rb,\ \ \l_{0}=t,\ \ \L_{-1}=\L'.$$
We have exact sequences of complexes
$$\xymatrix{0 \ar[r]& \Ct_{Z_p/\L_p,0}^{\cdot} \ar[r]& \Ct_{Z_p/\L_p,0}^{\cdot} \ar[r]&\Ct_{Z_{p+1}/\L_{p+1},0}^{\cdot} \ar[r]&0}$$
These exact sequences induce long exact sequences in cohomology
$$\xymatrix{\dots \ar[r]&H^k(\Ct_{Z_p/\L_p,0}^{\cdot}) \ar[r]&
H^k(\Ct_{Z_{p+1}/\L_{p+1},0}^{\cdot}) \ar[r]&H^{k+1}(\Ct_{Z_p/\L_p,0}^{\cdot}) \ar[r]&\dots.}$$
It is readily seen that the module $H^0(\Ct_{Z_p/\L_p,0}^{\cdot})$ can be identified with the module 
$\p_p^{-1}(\OM_{\L_p,0})$. Thus, in the exact sequence the map
$$\xymatrix{H^0(\Ct_{Z_p/\L_p,0}^{\cdot}) \ar[r]&
H^0(\Ct_{Z_{p+1}/\L_{p+1},0}^{\cdot})}$$
is surjective and therefore the exact sequence splits.
This shows that the induced maps
$$\frac{H^1(\Ct_{Z_p/\L_p,0}^{\cdot})}{\l_{p+1}H^1(\Ct_{Z_p/\L_p,0}^{\cdot})}
\to H^1(\Ct_{Z_{p+1}/\L_{p+1},0}^{\cdot})$$
are injective. Finally, the injectivity of those maps implies in turn that
the map $R$ is injective. This proves the assertion.\\
 
The coherence theorem for the direct image sheaves of the relative Lagrange
complex (Theorem \ref{T::coherent}, section \ref{S::coherent}) and the fact that
these sheaves are concentrated at the origin (Proposition \ref{P::concentrated}, section \ref{S::coherent}) imply that the $\OM_{\L',0}$-module $ H^1 (\Ct^{\cdot}_{Z'/\L',0}) $ is of
finite type. Therefore, using the previous assertion and the Nakayama lemma,
we get that the system of generators 
$$ \lb \bar \theta_F(\d_{\l_1}),\dots, \bar \theta_F(\d_{\l_k}) \rb $$
of the vector space $  H^1(\Ct^{\cdot}_{L,0}) $ lifts to a system of generators
$$ \lb \theta_G(\d_{\l_1}),\dots, \theta_G(\d_{\l_k}) \rb $$
of the $ \OM_{\L',0} $-module $ H^1 (\Ct^{\cdot}_{Z'/\L',0}) $.
This concludes the proof of Theorem \ref{T::versel}.

%%%%%%%%%%%%%%%%%%%%%%%%%%%%%%%%%%%%%%%%%%%%%%%%%%%%%%%%%%%%%%%
%%%%%%%%%%%%%%%%%%%%%%%%%%%%%%%%%%%%%%%%%%%%%%%%%%%%%%%%%%%%%%%
\section{The coherence theorem}
\label{S::coherent}
For the proof of the rigidity theorem to be complete, we need to prove 
the coherence of the direct image sheaves of the Lagrange complex and
the concentration at the origin of these sheaves. This is done in this section.
Our proof of the coherence follows that of Sevenheck and van Straten in the absolute case. It is similar to that of Greuel (\cite{Greuel}) for the de Rham complex, the idea of which goes back to Brieskorn (\cite{Br3}).
The main technical part to be proven is to adapt the Kiehl-Verdier theorem to
our situation. This is done in the appendix.
 
%%%%%%%%%%%%%%%%%%%%%%%%%%%%%%%%%%%%%%%%%%%%%%%%%%%%%%%%%%%%%%%%%
\subsection{Coherence of the direct image sheaves}
\label{SS::coherence}
% Let us fix some conventions. A topological manifold $ X $ of dimension $ d $
% will be called {\em compact with boundary} if $ X $ is a compact 
% topological space and if there exists a system of charts of $ X $,
% $(\p_i,U_i)  $ such that $ \p_i $ sends either $ U_i $ to $ \RM^d $ or 
% to a closed half space in $ \RM^d $.\\
We will say that a topological manifold with boundary $ X $ is {\em complex} if there exists a one 
parameter family of complex variety $ Y_t $ such that
\begin{enumerate}
\item $ Y_t \subseteq Y_{t'} $ for $ t'>t $,
\item $ X $ is equal to the closure of $ Y_0 $ in $ Y_t $ for any $t>0$.
\end{enumerate}
For a given holomorphic map
$$\p:Z \to \L  $$
we say that $ \p $ {\em extends} to $ \bar Z $ if $ \p $ is the
restriction of a deformation
$$ \p':Y_t\to \L $$
for any $ t>0 $. In such a case, we denote by 
$$\bar \p:\bar Z \to \L  $$
the map obtained by extending $ \p $ to the boundary of $ Z $.\\
Now we assume that $ \p $ is the deformation of a Lagrange variety.
Recall that the direct image image sheaves $ \RM^p \p_* \Ct_{Z/\L}^{\cdot} $ are
defined by the pre-sheaves
$$ U \mapsto \HM^p(\p^{-1}(U),\Ct_{Z/\L}^{\cdot}) .$$
\begin{theorem}
\label{T::coherent}
{ Let $\p:Z \to \L$ be a Lagrangian deformation such that
\begin{enumerate}
\item $ \p $ is pyramidal and extends to $ \bar Z $,
\item the fibres of $ \bar \p:\bar Z \to \bar \L $ are transverse to 
the boundary of $ Z $.
\end{enumerate}
Then, the direct image sheaves
$ \RM^p \p_* \Ct_{Z/\L}^{\cdot} $ are coherent
sheaves of $\OM_\L$ modules.}
\end{theorem}
For the case of a mapping to a point
$$ \p:L \to \lb 0 \rb $$
we recover the Sevenheck-van Straten finiteness theorem (\cite{VS}).
%%%%%%%%%%%%%%%%%%%%%%%%%%%%%
\subsection{Proof of the coherence theorem}
%\label{SS::vanStraten}
% To prove the theorem, we shall use a corollary of the Kiehl-Verdier 
% theorem which imitates a result of van Straten (\cite{vanStraten}).
% The proof of this corollary (that we state here as a theorem)
% is given in the appendix.\\
According to Theorem \ref{T::vanStraten} which is proved in the 
appendix, it sufficies to construct a {\em transversal datum $(\theta_i,U_i)$} to
$ (\p,\Ct_{Z/\L}^\cdot) $ i.e. a covering $ (U_i) $ of $ L $ together
with a collection of $ C^\infty $ vector fields $(\theta_i)$ such that
$\theta_i$ is transverse to the boundary of $ U_i $
and  the cohomology of the sheaves $ \Ct_{Z/\L}^\cdot $ are 
constant along the integral lines of $ \theta_i $.\\
The construction of such a field is local, it is 
therefore sufficient to consider the case where
\begin{enumerate}
\item $ Z=U $ is
a subvariety of some relatively compact open Stein space $ U  \subset (M \times \L) $,
\item $ \bar Z \subset (M \times \L)$ is compact,
\item the fibres of
$ \bar \p:\bar Z \to \bar \L $ intersect the 
boundary of $ \bar Z $ transversally.
\end{enumerate}
That $ \p $ is pyramidal implies that for any point $ x \in \d \bar Z $
there exists a holomorphic function $ h \in \It_x $ such that the
Hamilton vector field $ X_h $ of $ h $ is transverse to $ \d \bar Z $ and 
points towards $ Z $. Moreover, Proposition \ref{P::produit} 
(section \ref{SS::pyramidal}) implies that the sheaves $ \Ht^p(\Ct^\cdot_{Z/\L}) $ are constant
along the integral lines of $ X_h $. Consequently, we can take a covering
$ (V_k)  $ of the boundary of $ Z $ in $ M \times \L $  such that
there exists a transversal Hamiltonian vector field $ X_k   $
in $ V_k $.\\ 
 Let $ (\psi_k,V_k) $ be a partition of the unity and define the $ C^\infty $  vector-field
 $$ \theta=\sum_{k=1}^s\psi_k X_{k} .$$
This vector field is transversal
to $ (\p,\Ct_{Z/\L}^\cdot) $. 
This proves the theorem.
%%%%%%%%%%%%%%%%%%%%
\subsection{Concentration and freeness of the higher direct image sheaves}
%Another important algebraic property of this direct images is that in a small 
%neighbourhood of a singular point all the information is contained at 
%the singular point. This is an important property for explicit computations
%when dealing with 
%germs of deformations.\\
Consider a pyramidal Lagrangian deformation germ
$$ F:(\CM^k \times \CM^{2n},0) \to (\CM^n,0). $$ 
The deformation $\p:Z \to \L$ is called a {\em standard representative} of  
$F $ if  $ \p $ satisfies the assumptions of Theorem 
\ref{T::coherent}.
The existence of such a representative follows from the 
definition of a pyramidal deformation.
\begin{proposition}
\label{P::concentrated}
{Let $ \p:Z \to (\L,0) $ be a standard representative of the germ of a Lagrangian deformation of a pyramidal Lagrangian variety $ L $.
Then, there is a canonical isomorphism
$$ (\RM^p \p_* \Ct^\cdot_{Z/\L})_0 \approx H^p(\Ct^\cdot_{Z/\L,0}). $$}
\end{proposition}
The proof of the proposition is a straightforward variant of the corresponding
result due to van Straten for hypersurface singularities (\cite{VS}).
Consequently, it will not be given.
 
\newpage
% %%%%%%%%%%%%%%%%%%%%%%%%%%%%%%%%%%%%%%%%%%%%%%%%%%%%%%%%%%%%%
 \section{Appendix:The shrinking theorem}
 \label{S::KV}
 %%%%%%%%%%%%%%%%%%%%%%%%%%%%%%%%%%%%%%
\subsection{The Kiehl-Verdier theorem}
%Recall that a map between two complexes is called a {\em 
%quasi-isomorphism}
%if it is an isomorphism of the cohomology spaces.
\begin{theorem}
[\cite{Verdier},\cite{Douady}]
\label{T::Verdier}
{Let $A=(A_t)_{t \in [0,1]}$ be a nuclear chain of Fr\'echet
algebras. Let $E^{\cdot},F^{\cdot}$ be finite complexes of 
$A_0$-modules such that $A$ is
transverse to the $E^k$'s and the $F^k$'s.
Assume that there exists a morphism of complexes $f:E^{\cdot} \to 
F^{\cdot}$ such that
\begin{enumerate}
\item $f$ is a quasi-isomorphism,
\item the maps $f^k:E^k \to F^k$ are $A_0$-subnuclear
\end{enumerate}
Then there exists a complex  $L^{\cdot}$ of free $A_1$ modules of 
finite type
and a quasi-isomorphism $h:L^{\cdot} \to A_1 \tensor_{A_0} 
E^{\cdot}$.}
\end{theorem}
The definitions involved in the theorem (transversality, sub-nuclearity)
will not be given, we refer to \cite{DouadyR.}, \cite{Hubbard} and
\cite{Douady}. Whatever these definitions are is in fact unimportant 
to apply the theorem in the complex holomorphic situation.
% namely one
%has the following main corollary of the theorem.\\
%The space $\Ft(U)$ of global sections  of the sheaf
%$\Ft$ over a subset $U \subset (X \times S)$
%can be endowed with a Fr\'echet space 
%structure (\cite{Hubbard}).\\
%Let $s=(s_1,\dots,s_n), r=(r_1,\dots,r_n)$ be two sequences of real 
%positive numbers satisfying $s_j<r_j$. Denote by $D_t, t 
%\in [0,1]$  polydisks of radius $r-ts$ centered at the same point such that
%$D_0 \subset S$.\\
%
%Let $Y \bar{\subset} X$ be a relatively compact open Stein subset
%(recall that relatively compact means that the closure $\bar{Y}$ is  compact).
%\begin{proposition}
%\label{P::bbaki}
% {For a coherent analytic sheaf  $\Ft$ on $X \times S$, the following 
%properties
%hold
%\begin{enumerate} 
%\item The chain $(\OM_S(D_t))$ is a nuclear chain of Fr\'echet algebras
% transverse to $\Ft(D_0 \times U)$
%{\rm (\cite{Douady}, section 3, proposition 2 and section 6, example)},
%\item The restriction map  $r:\Ft(D_0 \times X ) \to \Ft(D_0 \times 
%Y)$
%is $ \OM_S(D_0) $-subnuclear
%{\rm (\cite{Douady}, section 4, proposition 4)},
%\item For any Stein open subset $U' \subset U$, there is a canonical
%isomorphism
%$ \OM_S(U') \tensor_{\OM_S(U)} \Ft(U \times X) \approx \Ft(U' \times X)$
%{\rm (\cite{Hubbard}, section 4, corollary 1)}.
%ù\end{enumerate}}
%\end{proposition}
Let $\Kt^{\cdot}$ be a complex of coherent sheaves on a complex
variety $S \times X$. Let $s=(s_1,\dots,s_n), r=(r_1,\dots,r_n)$ be
two sequences of real positive numbers satisfying $s_j<r_j$. Denote by $D_t, t 
\in [0,1]$  polydisks of radius $r-ts$ centered at the same point such that
$D_0 \subset S$.
 
Let  $Y \subset X$ be  relatively compact subset i.e. a subset for 
which the closure $\bar{Y}$ of $ Y $ in $ X $ is  compact.
By applying the Kiehl-Verdier theorem  with $A=(\OM_S(D_t))$, $E^{\cdot}=\Kt^{\cdot}(D \times X)$
and  $F^{\cdot}= \Kt^{\cdot}(D \times Y)$, we get the following 
complex analytic version of the Kiehl-Verdier theorem.
\begin{theorem}
\label{C::KV}
{Assume that the restriction map $r:E^{\cdot} \to 
F^{\cdot}$
is a quasi-isomorphism then there exists a complex $L^{\cdot}$ of free
$\OM_S(D_1)$-modules of finite type and a quasi-isomorphism
$h^{\cdot}:L^{\cdot} \to \Kt^{\cdot}(D_1 \times X)$.}
\end{theorem}
%This corollary is, of course, easier to handle than the general Kiehl-Verdier
%theorem.
%%%%%%%%%%%%%%%%%%%%%%%%%%%%%%%%%%%%%%%%
\subsection{The shrinking theorem}
Let $f:X \to S $ be a morphism between complex spaces and $\Kt^\cdot$
a complex of coherent sheaves on $X$
with a differential which is $ f^{-1}\OM_S $-linear.
\begin{definition}{\rm
The pair $(f,\Kt^\cdot)$ possess the {\em shrinking property}
if there exists a relatively compact subvariety $Y  \bar{\subset} X$
and open Stein coverings of $X$ and $Y$ 
$\underbar U=(U_i)$, $\underbar V=(V_i)$, where $V_i$ is a relatively compact subset of $ U_i$ and such that the restriction maps
$$C^p(\underbar U \cap f^{-1}(D),\Kt^q) \to C^p(\underbar V \cap f^{-1}(D),\Kt^q)$$
induce an isomorphism in hypercohomology
$$r_*:\RM^\cdot f_* \Kt^\cdot \to \RM^\cdot (f_{\mid Y })_*\Kt^\cdot.$$}
\end{definition}

%{\em Example.}{ Let $X$ be a compact manifold and let $\Kt^\cdot$
%be a complex consisting of one sheaf, say $\Ft$. Then, taking $Y=X$,
%we get that the Leray theorem implies $ (f,\Kt^\cdot)$ satisfies the
%shrinking property.}\\

\begin{theorem}
\label{T::shrinking}{If the map $f:X \to S$ satisfies the shrinking
property for a complex of sheaves $\Kt^\cdot$ then, the direct image 
sheaves $ \RM^\cdot f_*\Kt^\cdot $ are coherent.}
\end{theorem}

%%%%%%%%%%%%%%%%%%%%%%%%%%%%%%%%%5
\subsection{Proof of the shrinking theorem}
\ \\
Denote by $D \subset S$ a Stein open neighbourhood of $s \in S$.\\
Consider the complexes of Fr\'echet spaces
$$\xymatrix{E^{\cdot}:0 \ar[r]& E^0 \ar[r] \ar[d]^{r_0}& \dots 
\ar[r] &E^k \ar[r] \ar[d]^{r_k}& 0\\
F^{\cdot}:0 \ar[r]& F^0 \ar[r]& \dots \ar[r] &F^k \ar[r]& 
0}$$
where
\begin{enumerate}
\item $E^i=\bigoplus_{p+q=i}C^p(\underbar U \cap f^{-1}(D),\Kt^q)$,
\item $F^i=\bigoplus_{p+q=i}C^p(\underbar V \cap f^{-1}(D),\Kt^q)$, 
\item the differential in $E^{\cdot}$ and $F^{\cdot}$ equals
$d \pm \delta$, where $d$ is the \^Cech derivative and $\dt$ is the derivative
of the complex $ \Kt^{\cdot} $,
\item the map $r^\cdot$ is the restriction map.
\end{enumerate}
The cohomology vector spaces of $E^{\cdot}$ and 
$F^{\cdot}$ are the hypercohomology vector spaces
$\HM^\cdot(f^{-1}(D),\Kt^\cdot)$ and  $\HM^\cdot(f^{-1}(D) \cap Y,\Kt^\cdot)$.
Thus, by assumption the map $r$ is a quasi-isomorphism.\\
The complex $E^i$ is obtained as the space of sections over $X \times D$
of  a complex of sheaves in $X \times S$.
Namely, denote by $\Et^{p,q}$ the coherent
sheaf on $X \times S$ defined by the pre-sheaf
$$U \times U' \to C^p(U \cap f^{-1}(U')  , \Kt^q).$$
Put $\Et^i=\bigoplus_{p+q=i}\Et^{pq}$, then
$\Et^i(X \times D)=E^i$. Thus,
the complex version of the Kiehl-Verdier theorem applies (Theorem \ref{C::KV}).
Denoting by $D_0=D$ and $D_t$ polydisks like in the preceding 
subsection, we get that there exists
free  $\OM_S(D_1)$-modules of finite type $L^p$ and a quasi-isomorphism
$$h^{\cdot}:L^{\cdot} \to \Et^{\cdot}(X \times D_1) .$$
We sheafify this result as follows.\\
Denote by $k_p$ the rank of $L^p$ so that
$L^p \approx (\OM_S(D_1))^{k_p}$ and let $\Lt^p$ be the sheaf $ \OM_S^{k_p}$. Since $D_1$ is a Stein open subset,
the differential of the complex $L^{\cdot}$
induces a differential of the complex  $\Lt^{\cdot}$.\\
We assert that the quasi-isomorphism $h^\cdot$ extends to a 
quasi-isomorphism of sheaf complexes
$$\t h^\cdot:\Lt^\cdot \to \Et^\cdot(X \times -).$$
Indeed, let $U \subset D_1$ be a Stein open subset.
The map
$$1 \otimes h^p: \OM_S(U) \tensor_{\OM_S(D_1)}L^p \to  \OM_S(U) 
\tensor_{\OM_S(D_1)}
\Et^p(D_1 \times X)$$
is a quasi-isomorphism provided that $h^p$ is a quasi-isomorphism
(\cite{Douady}, corollary to Proposition 3 ).\\
Moreover, we have canonical identifications (\cite{Hubbard}, section 4, Corollary 1)
$$ \OM_S(U) \tensor_{\OM_S(D_1)}L^p=\Lt^p(U)$$
and 
$$ \OM_S(U) \tensor_{\OM_S(D_1)}\Et^p(X \times D_1)=\Et^p(X \times U).$$
This proves the assertion.
We get that
$$\Ht^p(\Lt^\cdot)(U) \approx \Ht^p(\Et)(X \times U)=\HM^p(f^{-1}(U),\Kt^\cdot)=(\RM^pf_*\Kt^\cdot)(U)$$
Thus, the direct image sheaves $\RM^p f_*\Kt^\cdot$ are coherent.
This concludes the proof of the shrinking theorem.
%%%%%%%%%%%%%%%%%%%%%%%%%%%%%
\subsection{A van Straten type criterion}
\label{SS::vanStraten}
We adapt the van Straten criterion for the coherence (\cite{vanStraten}) to more general 
situations.
% To prove the theorem, we shall use a corollary of the Kiehl-Verdier 
% theorem which imitates a result of van Straten (\cite{vanStraten}).
% The proof of this corollary (that we state here as a theorem)
% is given in the appendix.\\

Let $ \Kt^{\cdot} $ be a complex of coherent sheaves (bounded from 
the right) on a variety $ X $ and let $ f:X \to S $ be a holomorphic 
map. We assume that the differential of $ \Kt^{\cdot} $ is $f^{-1}\OM_S$
linear.
\begin{definition}[\cite{vanStraten}]{\rm A {\em transversal vector field}
to $ (f,\Kt^{\cdot}) $ in a Stein open subset $U \bar{\subset} X$
is a $ C^\infty $ vector field $ \theta $ defined in a 
neighbourhood of $\d \bar U $ in $ X $ satisfying the 
following 
conditions:
\begin{enumerate}
\item $ \theta $ is transversal to the boundary of $U$,
\item the integral curves of $ \theta $ are contained in the fibres of $ 
f_{\mid U} $,
\item the restriction of the cohomology sheaves $ \Ht^p(\Kt^{\cdot}) $
to the integral curves of $ \theta $ are constant sheaves.
\end{enumerate}}
\end{definition}
{\em Remark.}{ The last condition means that any cohomology class can be locally represented by a
set of holomorphic functions which are constant along the integral lines of $\theta$.}
\begin{definition}{\rm A {\em transversal datum} to $ (f,\Kt^{\cdot})$ is a
set $(\theta_i,U_i)$, $ i=1,\dots,n $, where $U_i$ is a Stein covering of $X$
and $\theta_i$ is a transversal vector field to
$(f,\Kt^{\cdot})$ in $U_i$.}
\end{definition}
\begin{theorem}
\label{T::vanStraten}{The existence of a transversal datum to $ 
(f,\Kt^{\cdot}) $ implies the existence of a shrinking. In 
particular, it implies that the hypercohomology sheaves $ \RM^p f_*\Kt^{\cdot} $ are coherent. }
\end{theorem}
The proof of this theorem repeats that of Theorem 1 in \cite{vanStraten}.

\bibliographystyle{amsplain}
\bibliography{master}
\end{document}